\documentclass[review]{article}
\usepackage[utf8]{inputenc}
\usepackage{amsmath}
\usepackage{amsfonts}
\usepackage{amssymb}
\usepackage{graphicx}
\usepackage{tabularx}
\usepackage{xstring}

\newtheorem{thm}{Theorem}

\newtheorem{claim}[thm]{Claim}

\newcommand{\lcm}{\textrm{lcm}}

\newcommand{\twoclass}[1]{\IfStrEq{#1}{}{}{\| #1\|_2}}

\newcommand{\gtwo}[2]{\twoclass{#1} > \twoclass{#2}}
\newcommand{\ltwo}[2]{\twoclass{#1} < \twoclass{#2}}
\newcommand{\eqtwo}[2]{\twoclass{#1} = \twoclass{#2}}
\newcommand{\netwo}[2]{\twoclass{#1} \neq \twoclass{#2}}
\newcommand{\getwo}[2]{\twoclass{#1} \ge \twoclass{#2}}

\begin{document}

\author{Robert Hochberg and Paul Phillips\\
Department of Mathematics\\
University of Dallas\\
1845 E. Northgate Drive\\
Irving, TX  75062\\
Phone: 972-721-5175, Fax: 972-721-4148\\
\{hochberg, phillips\}@udallas.edu}

\title{Discrepancy One among Homogeneous Arithmetic Progressions}

\maketitle

\begin{abstract}
We investigate a restriction of Paul Erd\H os' well-known problem from 1936 on the discrepancy of homogeneous arithmetic progressions. We restrict our attention to a finite set $S$ of homogeneous arithmetic progressions, and ask when the discrepancy with respect to this set is exactly 1. We answer this question when $S$ has size four or less, and prove that the problem for general $S$ is NP-hard, even for discrepancy 1.
\end{abstract}


\section{Introduction}

We are motivated by the Erd\H{o}s discrepancy problem \cite{erdos}. Given any sequence\\ $(x_n)\in \lbrace \pm 1\rbrace^\omega$, let us define
$d(s, k) = \sum_{i=1}^k x_{is}$. The problem asks whether for any such sequence $(x_n)$ and constant $C$ one may find $s$ and $k$ such that $\left |d(s, k)\right |\geq C$. We recall that a homogeneous arithmetic progression is one whose first term is equal to its common difference.  By defining the {\em discrepancy} of the homogeneous arithmetic progression starting at $s$ to be $\limsup\left\vert d(s, k) \right\vert$ as $k\to\infty$, we can recast the Erd\H{o}s problem in terms of a coloring of the positive integers with just two colors.  To wit, an answer of `no' would mean there is a way to color $\mathbb{Z}^+$ so that there is a universal upper bound for the discrepancy of every homogeneous arithmetic progression.  

This question goes back to the 1930s and has remained unsolved to this day. [Note that between submitting this paper for publication in July of 2015, and posting it to the archive in January of 2016, Terence Tao has announced a proof of the problem. See http://arxiv.org/abs/1509.05363.]
Mathias \cite{mathias} proved that if $d(s, k)\leq 1$ for all $s,k\geq 1$, then $d(s, k)$ has no lower bound. Konev et al. \cite{sat} showed that the interval $[1, 1160]$ could be 2-colored with discrepancy two, but any coloring of $[1, 1161]$ must contain a homogeneous arithmetic progression with discrepancy three.

In this paper we restrict the question in two ways: We first fix $C=1$ and then consider only {\em some} of the homogeneous arithmetic progressions. More precisely, given a ``skip set'' $S = \lbrace s_1, s_2, \ldots, s_n \rbrace$ of positive integers, what conditions on $S$ will guarantee that for any sequence  $(x_n) \in \lbrace \pm 1\rbrace^\omega$ there is at least one $s\in S$ and some positive integer $k$ so that $\left\vert d(s, k)\right\vert > 1$?  In such a case, we shall say that skip set $S$ \emph{forces discrepancy two}.

In section~\ref{Paths} we provide definitions for certain key concepts, in particular, the Divisibility and Parity Conditions and the Cycle Equation.  We also prove a Master Lemma and derive many of  the basic tools for our theory from it.    This lemma will then be used in section~\ref{OddCycles} to completely characterize skip sets of size four or less that force discrepancy two.  Finally, in section~\ref{Complexity} we will show that the general question of characterizing skip sets that force discrepancy two is NP-hard.  

\subsection{Reformulation of the Problem}
\begin{center}
\begin{figure}[h]
\centering
\includegraphics[scale=0.5]{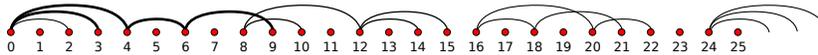}
\caption{The skip graph $G(2, 3, 4)$}\label{FIGg234}
\end{figure}
\end{center}

To investigate this question, we shall pass from a consideration of the sequences $(x_n)$ to an underlying graph $G(S)$ on the natural numbers determined by the skip set $S$.  The vertices of $G(S)$ are the natural numbers $\mathbb{N}=\{0, 1, 2, \ldots\}$.  For each $s \in S$ the graph has an edge between every even multiple of $s$ and the subsequent odd multiple of $s$. An edge arising from skip size $a$ is called an $a$-{\em arc}. $G(\{2, 3, 4\})$ is shown in Figure~\ref{FIGg234}. Note that this graph is periodic, with period $24 = 2\cdot\lcm{(2, 3, 4)}$.

To have a 2-coloring with discrepancy one for all the skips in the skip set, the vertices of each edge need to have opposite colors. Thus a skip set forces discrepancy two exactly when it produces an odd cycle in the underlying graph.
(Note: There is a small technical difference between this formulation and the Erd\H{o}s formulation given above in that our underlying graph starts at zero, while there is no $x_0$ term in the $\pm 1$ sequences.  As we are dealing with finite skip sets, we can define $L = \lcm(S)$.  It is then clear that the underlying graph for $[0,2L-1]$ and the corresponding graph for the actual Erd\H{o}s problem on $[1,2L]$ are just mirror images of each other.)
We wish, therefore, to determine necessary and sufficient conditions on a skip set $S$ so that $G(S)$ contains an odd cycle.  

\subsection{Notation}

We define equivalence classes 
based on how many factors of 2 a given positive integer has.  We denote the equivalence class of $a$ by $\twoclass{a}$.  For two positive integers $a$ and $b$ the notation $\gtwo{a}{b}$ indicates that $a$ has more factors of 2 than $b$ does, while $\eqtwo{a}{b}$ indicates the two numbers are in the same equivalence class. All odd numbers lie in the lowest equivalence class.

We use brackets, as in $[abcd]$, to indicate that a certain sequence of skip sizes occurs consecutively on a path in the underlying graph, 
without specifying whether the skips are to be taken in the positive or negative direction. To specify direction we will use pluses and minuses, so that $[-3+4+2+3]$ would indicate the bold path shown in the preceding figure, as would $[-3-2-4+3]$. When a pattern is realized on the number line, the values taken on by the endpoints of the arcs are called {\em terms} of the pattern. If a term falls at the right end of two arcs or the left end of two arcs, then we call that term a path {\em reversal}. A pattern of $k$ skips would have $k+1$ terms, and a cycle of $k$ skips would have $k$ terms. If a pattern (whether given with signs, without signs, or abstractly, as in $[abcab]$) can be found in some skip graph without repeating any terms or arcs, then we call that pattern {\em realizable}. If the pattern can be found possibly repeating some terms and/or arcs, then we call that pattern {\em weakly realizable}.  For example, $[2\ 1\ 3\ 4]$ is weakly realizable, but not realizable, and $[5\ 1\ 10]$ is not even weakly realizable. Patterns that are not realizable are called {\em forbidden}.

\section{Paths}\label{Paths}

This section concerns paths that occur in our graph. All of these claims apply to cycles as well.

Given the pattern $[\epsilon_1 a_1 + \epsilon_2 a_2 + \cdots + \epsilon_n a_n]$ where each $\epsilon_i \in \{\pm 1\}$, we define the \textbf{Intermediate Path} for this pattern to be $$I = [\epsilon_2 a_2  + \cdots + \epsilon_{n-1} a_{n-1}].$$  By a slight abuse of notation, which we trust context will keep unambiguous, we will also denote the \textbf{sum} of this intermediate path by $$I = \sum_{i=2}^{n-1} \epsilon_i a_i.$$ 
 
The \textbf{Divisibility Condition} is satisfied for this pattern if $I$ is a multiple of $\gcd(a_1, a_n)$.

Given that the Divisibility Condition holds, we define the \textbf{Parity Condition}  to be that $I$ is an  even multiple of $\gcd(a_1,a_n)$ if and only if
\begin{eqnarray*}	\ltwo {a_1}{a_n} &\textrm{implies}& \epsilon_1 = -1\\
 	\gtwo{a_1}{a_n}  &\textrm{implies}&  \epsilon_n = +1\\
 	\eqtwo{a_1}{a_n}  &\textrm{implies}& \epsilon_1 = - \epsilon_n
\end{eqnarray*}
 
Note that all possible relations between the equivalence classes of $a_1$ and $a_n$ are covered and that in each case there is no constraint on the sign of the path element in the higher equivalence class.

We say that a path satisfies the \textbf{Basic Parity Test}, to be abbreviated BPT, if every subpath consisting of two consecutive path elements satisfies the Parity Condition.  Claims~\ref{plusplus}, \ref{agtb}, \ref{plusminus}, and \ref{aeqb} below summarize the requirements for BPT. Later we will see that under certain conditions BPT will be sufficient to assert that a pattern is (weakly) realizable.

Our first claim gives us an exact description for when a path is realized in the underlying graph.

\begin{claim}[Master Lemma]\label{masterlemma}  
A path $[\epsilon_1 a_1 + \epsilon_2 a_2 + \cdots + \epsilon_n a_n]$ is weakly realizable if and only if 
each subpath satisfies the Divisibility and Parity conditions given above.
\end{claim}

\paragraph{Proof:} Pick any subpath and represent it as $[\epsilon_i a_i + \cdots + \epsilon_j a_j]$ for some $i<j$.  Denote the finishing term for the $\epsilon_i a_i$-arc by $F$ and the beginning term for the $\epsilon_j a_j$-arc by $B$.  Note that $F$ and $B$ are both multiples of $g=\gcd(a_i,a_j)$.  

\emph{Part I, Necessity:}
The sum $I$ of the intermediate path must be a multiple of $g$ as $I=B-F.$ We now consider the three possible cases for the respective equivalence classes of $a_i$ and $a_j$.

\emph{Case $\gtwo{a_i}{a_j}$}:  

The gcd $g$ has only as many factors of 2 as the representative of the lesser equivalence class, so $\eqtwo{g}{a_j}$.  Because $F$ is a multiple of $a_i$, we also have $\getwo{F}{a_i}\gtwo{}{g}$.  As $\epsilon_j = +1$ means $B$ is the left term of the $a_j$-arc, $B$ is an even multiple of $a_j$, and so $\epsilon_j=+1$ if and only if $\gtwo{B}{a_j}\eqtwo{}g$.  Thus $\epsilon_j = +1$ if and only if $\eqtwo{I}{B-F} \gtwo{}{g}$, which is to say, exactly when $I$ is an even multiple of $g$.

\emph{Case $\ltwo{a_i}{a_j}$}:

This case follows by applying the previous case to $ [-\epsilon_j a_j -  \cdots  - \epsilon_i a_i] $ which simply traces the path backwards.

\emph{Case $\eqtwo{a_i}{a_j}$}:

Note that $\eqtwo{g}{a_i}\eqtwo{}{a_j}$.  Here $\epsilon_i = -\epsilon_j$ if and only if $B$ and $F$ are both left vertices or both right vertices for their respective arcs.  Thus opposite $\epsilon$ values occur if and only if $B$ and $F$ are both be even multiples of $g$ or are both odd multiples of $g$, which is exactly when $I = B-F$ will be an even multiple of $g$.  

\emph{Part II, Sufficiency:}

We proceed by induction.  

We begin by noting that a path $[\epsilon_1 a_1 + \cdots + \epsilon_n a_n]$ is weakly realizable if and only if the following system of modular equations is solvable for $T$.  Here $T$ is the starting vertex for the given path.  The $i^\textrm{th}$ equation checks that the starting vertex for the $\epsilon_i a_i$-arc is an even multiple of $a_i$ when $\epsilon_i=+1$ and an odd multiple otherwise.
	\begin{align*}	
		T &\equiv \frac{1-\epsilon_1}{2} a_1  &\mod 2a_1 \\
		T &\equiv\frac{1-\epsilon_2}{2} a_2 - \epsilon_1 a_1  & \mod 2a_2 \\
		&\vdots&\\
		T &\equiv \frac{1-\epsilon_n}{2} a_n - (\epsilon_1 a_1 + \cdots + \epsilon_{n-1} a_{n-1}) &\mod 2a_n 		
	\end{align*}
	
By the generalized Chinese Remainder Theorem, a system of modular equations of this sort is solvable for $T$ if and only if the pairwise differences of the residues of $T$ are divisible by the greatest common divisor of the two respective moduli.

If $n=1$, then there is obviously no problem solving for $T$.  So suppose for induction that there is an integer $k\ge 1$ such that any (signed) path with $k$ path elements whose every subpath satisfies the Divisibility and Parity conditions has a corresponding set of $k$ modular equations that are simultaneously solvable for $T$.  We now consider a (signed) path of length
 $n=k+1$ satisfying the Divisibility and Parity conditions for all of its subpaths.  

The first $k$ steps of the path satisfy the conditions, hence the first $k$ equations are simultaneously solvable for $T$.    Likewise, the last $k$ steps of the path satisfy the conditions, so again by the induction hypothesis, the last $k$ equations are simultaneously solvable for $T + \epsilon_1 a_1$, hence solvable for $T$.  Combining these two observations lets us conclude that every pair of residues for $T$, except possibly the first and the last, has a difference divisible by the greatest common divisor of the two corresponding moduli.

Thus it only remains to check that the first and last equations are compatible.  That is, we need to check that 
\begin{eqnarray*}
  \frac{1+\epsilon_1}{2} a_1 - \frac{1-\epsilon_n}{2} a_n + (\epsilon_2 a_2 + \cdots + \epsilon_{n-1} a_{n-1}) &\equiv& 0 \qquad  \mod\gcd(2a_1,2a_n)
\end{eqnarray*}

However, by the stated conditions, the term in parentheses, being the sum of the intervening path, is a multiple of the $\gcd(a_1, a_n)$.  As long as it has the same parity as 
$\frac{1+\epsilon_1}{2} a_1 - \frac{1-\epsilon_n}{2} a_n $ does with respect to $\gcd(a_1, a_n)$, then the equation is true.  This is exactly what the Parity condition supplies.

The induction being complete, we have proved sufficiency.\hfill$\Box$

\paragraph{Corollaries from the Master Lemma.}  The following corollaries are presented as they are useful and fairly easily remembered.  Where they are simply special cases of the Master Lemma, proofs are omitted. Many rely on the observation that if a pattern is not weakly realizable, then it is not realizable, and therefore forbidden.

\begin{claim}[Basic Validity]\label{aa}
The pattern $[aa]$ is forbidden in any path.
\end{claim}

\paragraph{Proof:}  Because $\eqtwo aa$, the two $a$-arcs have opposite signs, thus the path element $[aa]$ repeats an arc, which we disallow.
\hfill$\Box$

\begin{claim}\label{aba}
The pattern $[aba]$ in a path implies $a\mid b$.
\end{claim}

\begin{claim}\label{abab}
The pattern $[abab]$ is forbidden in any path.
\end{claim}

\noindent\textbf{Proof:} By the previous claim we would have $a\mid b$ and $b\mid a$, implying $a=b$.
\hfill$\Box$

\begin{claim}\label{abageneral}
Generalizing Claim~\ref{aba}, if the pattern $[ax_1x_2\ldots x_kb]$ appears in any path, then $\gcd(a,b) \mid x_1\pm x_2\pm\ldots \pm x_k$, for some appropriate choice of additions and subtractions.  And in the special case where $a=b$ we have $a \mid x_1\pm x_2\pm\ldots \pm x_k$. 
\end{claim}

\begin{claim}\label{plusplus}
If the pattern $[+a+b]$ or $[-b-a]$ appears in any path, then $\gtwo{a}{b}$.

\end{claim}

\begin{claim}\label{agtb}
Given $\gtwo{a}{b}$ then the following are all the realizable path segments where $a$ and $b$ are neighbors:
$[\pm a+b]$ and $[-b \pm a]$.

\end{claim}

\begin{claim}\label{plusminus}
If the pattern $[+a-b]$ appears in any path, then $\eqtwo{a}{b}$.
\end{claim}

\begin{claim}\label{aeqb}
 Given $\eqtwo{a}{b}$ where $a$ and $b$ are consecutive path elements, then they must have opposite signs.
\end{claim}

\begin{claim}\label{oddblocks}
In any path the pattern $[a x_1 x_2 \ldots x_k b]$ where $a$ and $b$ (not necessarily distinct) are even and all of the $x_i$ are odd, we must have that $k$ is even. That is, odd skips must come in even-size blocks.
\end{claim}

\begin{claim}\label{scalePattern}
For any positive integer $d$, the pattern $[x_1x_2\ldots x_k]$ is (weakly) realizable if and only if the pattern $[(dx_1) (dx_2) \ldots (dx_k)]$ is (weakly) realizable.
\end{claim}

\noindent\textbf{Proof:} The conditions in the Master Lemma are not affected by adding or removing a common factor to/from every skip in a pattern.
\hfill$\Box$

\noindent\textbf{Other path results that will be used later:}

\begin{claim}\label{abca}
If $a>b$ and $a>c$, then the pattern $[abca]$ is forbidden.
\end{claim}

\noindent\textbf{Proof:} By Claim~\ref{abageneral} some combination $\pm b\pm c$ must give a multiple of $a$.  Because both $b$ and $c$ are smaller than $a$, we have $0 < \vert b-c \vert < a$ and $0 < b+c < 2a$, and thus $b+c=a$.  Both $b$ and $c$ must thus be given the same sign in the pattern.  Furthermore, if either instance of $a$ had the opposite sign then the pattern would contain a 3-cycle, and therefore repeat a term.  Thus, the pattern must be either $[+a+b+c+a]$ or $[-a-b-c-a]$. But by Claim~\ref{plusplus} applied three times, we would have $\gtwo{a}{a}$, a contradiction.\hfill$\Box$

\begin{claim}\label{aabbcc}
Any permutation of the pattern $[aabbcc]$ is forbidden.
\end{claim}

\noindent\textbf{Proof:} Without loss of generality we may assume that the pattern begins with $[ab]$. The third term may be $a$ or $c$.

{\em Case $[aba]$}: We must place $\{b, c, c\}$ to finish the pattern, without the $c$'s being together, and the only way to do that is with the pattern $[abacbc]$. But then we have $a\mid b$ and $c\mid b$, implying that $a<b$ and $c<b$. But since $[bacb]$ is part of that pattern, we have a contradiction with Claim~\ref{abca}.

{\em Case $[abc]$}: By Claim~\ref{aa} there are four sub-cases, which we treat separately.

{\em Subcase $[abcabc]$}: Suppose that $a > b$ and $a > c$. Then by considering the sub-pattern $[abca]$ we get a contradiction with Claim~\ref{abca}. The same is true whether $b$ or $c$ is the largest skip size.

{\em Subcase $[abcacb]$}: Here we have that $c\mid a$, and either $a\mid b+c$ or $a\mid b-c$, which in either case yields $c\mid b$. Thus $c$ is a divisor of each term, so by Claim~\ref{scalePattern} we may assume that $c=1$, giving the pattern $[ab1a1b]$. By Claim~\ref{oddblocks}, $b$ must be odd, and by Claims~\ref{abageneral} and~\ref{plusminus}, $a\mid b-1$, and therefore $a<b$. The pattern $[b1a1b]$ implies that the skips $\pm 1\pm a\pm 1$ must span the distance between two multiples of $b$, and since $a < b$, those three skips must all have the same sign. But then by applying Claim~\ref{plusplus} twice we get $\gtwo 1 1$, a contradiction. 

{\em Subcase $[abcbac]$}: This pattern is equivalent to tracing the previous pattern in reverse. Because we made no assumptions about the signs of the skips in that pattern, the proof is the same.

{\em Subcase $[abcbca]$}: This contains the pattern $[bcbc]$ which cannot occur by Claim~\ref{abab}.\hfill$\Box$

\begin{claim}[No Pairs Corollary]\label{nopairs}
No pattern can contain any permutation of one, two, or three pairs.  That is, any permutation of $[xx]$, $[xxyy]$ or $[xxyyzz]$ is forbidden.
\end{claim}

\noindent\textbf{Proof:}  This is a direct corollary of Claims~\ref{aa}, \ref{abab}, and~\ref{aabbcc}.\hfill$\Box$

On the other hand, we note that a permutation of four pairs is possible; for instance, the pattern $[-1-4+6+3-1+4+6+3-1-4]$ can be realized by making 5 the first term for the pattern.

\begin{claim}\label{bacabac}
The pattern $[bacabac]$ cannot occur in any path.
\end{claim}

\noindent\textbf{Proof:}  We invoke Claim~\ref{scalePattern} to reduce to the case where $a, b$ and $c$ are not all even. By Claim~\ref{aba} we have $a\mid b$ and $a\mid c$ so that $\getwo b a$ and $\getwo c a$, and thus $a$ is odd. Here $b$ must also be odd, otherwise the sub-pattern $[bacab]$ would have an odd number of odd arcs between two consecutive even arcs, contradicting Claim~\ref{oddblocks}. Similarly, $c$ must be odd. Since all skips are odd,  by Claim~\ref{aeqb} the signs of the skips must alternate.  Let us assume without loss of generality that the pattern is $[+b-a+c-a+b-a+c]$. Because $b$ and $c$ are multiples of $a$, we have $b, c \geq 2a$. By Claim~\ref{abageneral} we have $b\mid (c-2a)$ and $c\mid (b-2a)$, which cannot be true as then $b$ and $c$ would each be smaller than the other.\hfill$\Box$

The following claim shows that the Parity condition for paths of length 3 follows from the Divisibility condition and the Basic Parity Test. It will considerably shorten our proofs about 5-cycles in the next section.

\begin{claim}\label{SingleSeparatedNeighbors}
If the path pattern $[+\epsilon_a a +\epsilon_b b +\epsilon_c c]$ satisfies satisfies the Basic Parity Test, then the Divisibility condition for the pattern implies the Parity condition.
\end{claim}

\noindent\textbf{Proof:}  By tracing the path backwards if need be, we may assume without loss of generality that $\getwo a c$.  Thus for $g=\gcd(a,c)$ we have that $\eqtwo g c$.  Furthermore, by Divisibility  $b$ is a multiple of $g$, so we know that $\getwo b g \eqtwo{}c$.  This gives us two cases to investigate: either $\gtwo b c$ or $\eqtwo b c$.

\emph{Case $\gtwo bc$:}  In this case, by BPT, $\epsilon_c = +1$.  Further, $b$ is an even multiple of $g$, so to satisfy Parity, $\epsilon_a$ and $\epsilon_c$ must behave as if $a$ and $c$ were consecutive path elements.  If $\gtwo ac$, then we are done as $c$ has the correct sign.  Otherwise, $\eqtwo ac$, in which case because $a$ precedes $b$ in the path, BPT forces $\epsilon_a=-1$ and so $a$ and $c$ have opposite signs as required.

\emph{Case $\eqtwo bc$:}  The Parity condition puts a requirement on the sign of the path element in the lower equivalence class, depending on whether the intervening path sum is an even or odd multiple of the gcd. The intervening path for  $[+\epsilon_a a +\epsilon_b b]$ is 0 (an even multiple) and the intervening path for $[+\epsilon_a a +\epsilon_b b +\epsilon_c c]$ is $\epsilon_b b$, an odd multiple. Therefore $\epsilon_c$ should equal $-\epsilon_b$, which is exactly what BPT specifies when $b$ and $c$ are in the same equivalence class.\hfill$\Box$

We note that this does not extend to longer paths, as shown by the sequence $[+4+3-1+2]$ which does not satisfy the Parity condition, even though it and all its subpaths satisfy Divisibility and it satisfies BPT.

\subsection{An Application to Cycles}

The Master Lemma also gives conditions on the realizability of cycles.  In addition to the Divisibility and Parity conditions given there, we also need the \textbf{Cycle Equation} 
$$\sum \epsilon_i a_i = 0,$$ 
to hold for the proposed cycle.

Moreover, there are now two intervening paths between any two path elements.  For any $i<j$, we can recast the cycle 
$[+\epsilon_1 a_1 + \epsilon_2 a_2 + \cdots + \epsilon_n a_n]$ as 
$[+\epsilon_i a_i + I_{ij} + \epsilon_j a_j + K_{ij}].$  The $I$ corresponds to the original intervening path between the $i^{\textrm{th}}$ and $j^\textrm{th}$ path elements as given in the Master Lemma, and $K$ is the other path that arises now that the original pattern is considered to be a cycle.

Conveniently, the Cycle Equation guarantees that the Divisibility and Parity conditions hold for $[+\epsilon_i a_i + I+ \epsilon_j a_j]$ if and only if they hold for $[+\epsilon_j a_j + K+ \epsilon_i a_i]$, and thus when verifying the conditions of the Master Lemma on a cycle, we need only check whichever subpath is convenient.

\begin{claim}\label{BothPaths}
Suppose $I$ and $K$ are (signed) paths and the Cycle Equation holds for the pattern $[+\epsilon_a a +I+ \epsilon_b b+ K]$.  Then the Divisibility and Parity conditions hold for $[+\epsilon_a a +I+ \epsilon_b b]$ if and only if they hold for $[+ \epsilon_b b+ K+\epsilon_aa]$.  
\end{claim}

\noindent\textbf{Proof:}
From the Cycle Equation we have that $K = - (\epsilon_a a + I + \epsilon_b b)$.  Hence, if $g=\gcd(a,b)$ divides either $I$ or $K$, then it necessarily divides the other.

Recall that if $\eqtwo ab$ then $g$ is in the same equivalence class as $a$ and $b$, hence $a$ and $b$ are each odd multiples of $g$.  On the other hand, if $a$ and $b$ are in different equivalence classes, then one is an odd multiple of $g$ and the other is an even multiple.  Thus to check the Parity conditions, it suffices to observe that the Cycle Equation guarantees that $\eqtwo ab$ if and only if $I\equiv K\pmod{2g}$.   

Indeed, if $\eqtwo ab$, then because $I$ and $K$ are congruent modulo $2g$, the Parity conditions from the two intervening paths impose the same conditions on the relative values of $\epsilon_a$ and $\epsilon_b$.

On the other hand, if $\netwo ab$, then $I= -K\mod 2g$.  Because the path element from the lesser equivalence class is on opposite sides of these intervening paths, opposite aspects of the Parity condition are invoked for them.  However, as $I$ and $K$ have opposite parities with respect to $g$, these opposite aspects yield the same requirement.\hfill$\Box$

\section{Odd Cycles}\label{OddCycles}

We are interested in 2-colorings of the natural numbers that have discrepancy 1 with respect to a given set of skip sizes. Such a 2-coloring exists if and only if the underlying graph has no odd cycles. So we will now consider the question of skip sets that generate odd cycles. We leave it to the reader to show how Claim~\ref{aa} implies the following theorem.

\begin{thm}[Skip Sets of Sizes 1 and 2]
If $|S|\leq 2$ then $\mathbb{N}$ can be 2-colored with discrepancy 1, with respect to each of the homogeneous arithmetic progressions with skip sizes in $S$.
\end{thm}

\begin{claim}\label{lowestOddClass}
If $S$ is a reduced skip set with $|S|\ge 3$ and $S$ generates an odd cycle using each skip at least once, then there are at least two odd skips in $S$.
\end{claim}

\noindent\textbf{Proof:}  Because $S$ is reduced, at least one of the skips is odd, so that proving the claim consists of showing that at least one other skip size must be odd.  The rest of the skips cannot all be even, however, else the sole odd skip would need to be paired with itself to satisfy Claim~\ref{oddblocks}, an impossibility by Claim~\ref{aa}. \hfill$\Box$

\begin{claim}\label{largesttwice}
If $S$ is a skip set with $|S|\le 4$, then the largest skip size cannot appear more than once in any cycle.
\end{claim}

\noindent\textbf{Proof:} We begin by invoking Claim \ref{scalePattern} so that we can assume that $S$ is reduced. Suppose for contradiction that $a$ is the largest skip size and that it appears at least twice in the cycle. Consider the $a$-arcs of our cycle as they lie on the number line, and select two consecutive ones, $L$ and $R$, with $L$ coming before $R$ on the number line.  Note that $L$ and $R$ cannot share vertices but the vertices are all multiples of $a$.  Letting $J$ denote the interval (and its length) between these arcs we see that its length is a positive multiple of $a$, and so $J \geq a$.

These two arcs must be connected by two paths constructed from the skip sizes $b$, $c$ and $d$.  Because these skip sizes are all less than $a$, those two paths much have at least one term each in the interior of  interval $J$.  

To see that no term in this interval can be a path reversal, we note that a path reversal term would be a multiple of two different values, say $b$ and $c$, and yet both the path containing the reversal term and the other path will have to contain an arc that passes over this term in order to traverse $J$. As all $b$- and $c$-arcs will land on the reversal term or stay to one side of it, only a $d$-arc will suffice.  However, there is at most one $d$-arc that passes over this term, and only one of the two paths may contain it.  Thus any path reversal in $J$ would preclude at least one of the two paths from traversing the interval $J$.

Now let $x$ and $y$, with $x<y$, be two consecutive terms in $J$, one from each path connecting $L$ to $R$; that is, no other terms of the cycle lie between them on the number line. Because the terms are not path reversals, we can assume without loss of generality that $[+b+c]$ is the pattern of arcs meeting at term $x$. The $c$-arc must pass over $y$, leaving only $b$ and $d$ to meet at $y$.  At $y$ we can rule out $[+b +d]$ as that would require the $b$-arc to pass over $x$ which is impossible as $x$ is a multiple of $b$.  Thus $[+d+b]$ must be the pattern of arcs meeting at term $y$, and by Claim~\ref{plusplus}, $\gtwo{d}{}\gtwo{b}{c}$.  Because Claim~\ref{lowestOddClass} requires us to have two odd numbers, both $a$ and $c$ must be odd. 

We get our contradiction by considering the right end of the $c$-arc incident with $x$. Because $c$ is odd Claims~\ref{agtb} and~\ref{aeqb} show that the next arc must be odd and the corresponding term a path reversal.  The only other odd skip size is $a$ and so this gives that the pattern must be $[+b+c-a]$.  There are no path reversal terms in $J$, so this term must be to the right of interval $J$ as well as being the right endpoint of an $a$-arc.  However, because $c<a$ this reversal term must lie in the region between the two endpoints of the $a$-arc named $R$, where there are no multiples of $a$.  This is a contradiction.\hfill$\Box$

On the other hand, we note that the pattern 
$[+18+9-3+16+20+3-9-18+3-19-20]$ on skip set
$S=\{3, 9, 16, 18, 19, 20\}$ 
generates an 11-cycle starting at 360, and uses its largest skip size twice.

\begin{thm}[Skip Sets of Size 3]\label{skipsetsize3}
A skip set $\{a, b, c\}$ of size 3, with $a,b<c$, will force discrepancy two if and only if $a+b=c$ and $\netwo ab$.  Moreover, the underlying graph will contain 3-cycles, but no other odd cycles.
\end{thm}

\noindent\textbf{Proof:} 
Suppose first that the set $\{a,b,c\}$ forces discrepancy two. Then $G(\{a, b, c\})$ contains an odd cycle.  By Claim~\ref{largesttwice}, $c$ cannot appear twice in this cycle. And since the patterns $[aa]$ and $[abab]$ are forbidden in any path or cycle by the No Pairs Corollary (Claim~\ref{nopairs}), it is impossible to create an odd cycle of size greater than 3. Thus any odd cycle must be a triangle of the form $[+a+b-c]$ or $[+b+a-c]$. The Cycle Equation gives $a+b=c$, and by Claim~\ref{plusplus}, $\netwo ab$.

Now suppose that we have a set $\{a, b, c\}$ satisfying $a+b=c$ and $\netwo ab$. We assume without loss of generality that $\gtwo{a}{b}$ and consider the candidate cycle $[+a+b-c]$.  
The Cycle Equation, $a+b=c$, then shows that $\eqtwo bc$ and also that Divisibility holds for every path in the proposed cycle. 
By Claim~\ref{SingleSeparatedNeighbors}, because the largest separation of path elements is by a single other element, it only remains to check the BPT.  This follows immediately, however, by application of Claims~\ref{aeqb} and~\ref{agtb}. Finally we note that the odd cycle implies that $\{a, b, c\}$ forces discrepancy two.
\hfill$\Box$

\begin{claim}\label{fourskipsizes}
A graph generated by four skip sizes cannot contain an odd cycle of length greater than 7.
\end{claim}

\begin{center}
\begin{figure}
\centering
\includegraphics[scale=.2]{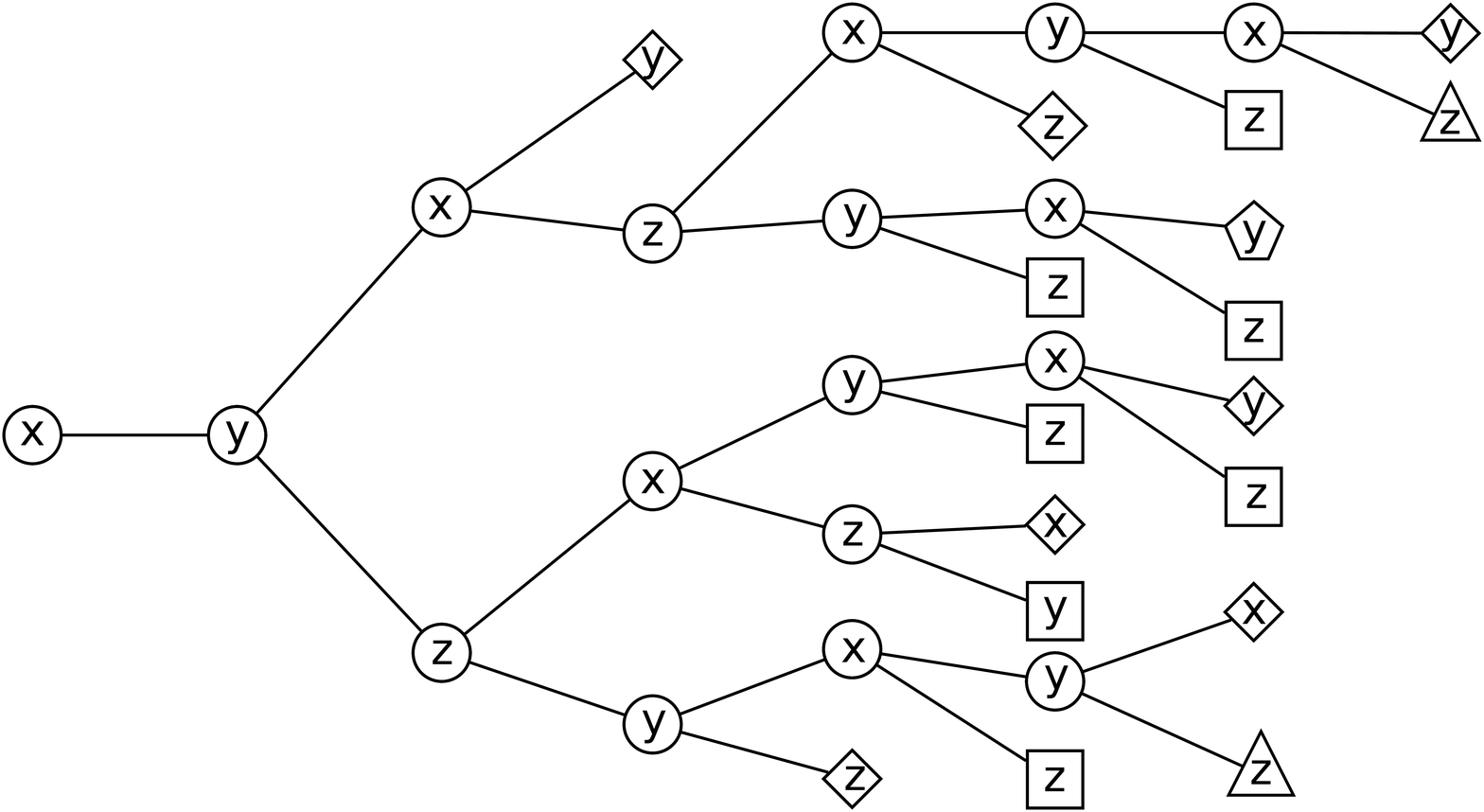}
\caption{Proof of Claim~\ref{fourskipsizes}}\label{treeproof}
\end{figure}
\end{center}

\noindent\textbf{Proof:} 
Let $a$ be the largest skip size. By Claim~\ref{largesttwice} this skip cannot occur again in the cycle.  We will now see that no path of length greater than 7 can be made from the three remaining skip sizes $\{x, y, z\}$. Suppose some such path begins $[xy]$. Figure~\ref{treeproof} shows all ways to extend this path subject to Claim~\ref{aa}. Nodes with a diamond violate Claim~\ref{abab}, nodes with a square violate Claim~\ref{aabbcc}, the pentagon node requires a path where $x=y$, as they must divide each other, and nodes with a triangle violate Claim~\ref{bacabac}, so that all paths of length 9 or longer beginning with $a$ are forbidden.
\hfill$\Box$

\begin{claim}\label{SevenCycle} A reduced skip set of size four generates a 7-cycle if and only if the elements can be labeled $\lbrace a, x, y, z\rbrace$ so that all the following conditions are met.
	\begin{itemize}
		\item  $z = 1$; 
		\item  $x$ and $y$ are odd; 
		\item  $a = 2x + y - 3$ (making $a$ even); 
		\item  $x \mid (a+1)$; and
		\item $\gcd(a,y) = 1$.
 	\end{itemize}
\end{claim}

\noindent\textbf{Proof:}
\emph{Part I, Necessity:} By Claim~\ref{aa} no element may appear four times in a 7-cycle, and by Claim~\ref{aabbcc} the cycle cannot contain three elements that each appear twice. Thus some element must occur three times, another twice, and the other two once each.  By spacing considerations and the fact that the pattern $[sts]$ forces $s\mid t$, we see that the triply repeated element is the gcd of the skip set, hence is equal to 1.  Let $z$ be this element and the first requirement is shown to be necessary.  

So without loss of generality, the pattern is $[a1x1y1x]$.  As each 1 in the cycle is required to be paired with another odd number in the path (Claim~\ref{oddblocks}), $x$ must be odd, as well as exactly one of $a$ or $y$.  


To see that $y$ cannot be even, note that if it were, then the cycle must be $[+a -1 + x -1 \pm y + 1 -x ]$ and thus we obtain $x$ divides $ -1  \pm y +1 = \pm y$.  But then $\gcd(x,y) = x$ and because of the path element $[y1x]$ we need $x \mid 1$, but $x \ne 1$.  

Thus $a$ is the only even number and we arrive at the second necessary requirement that $x$ and $y$ are odd.

As $a$ is even and the rest of the skip numbers are odd our signed pattern must conform to $[\pm a + 1 - x + 1 - y + 1 - x]$.  Yet starting with $-a$ is impossible, as then we would need $a = 3 - 2x -y$, which of necessity would be negative in order to make the pattern a cycle.
Thus, in fact, the signed pattern is $[+a + 1 - x + 1 - y + 1 - x]$, and the cycle condition will require the third necessary condition that $a=2x + y -3$.  Moreover, by the master lemma, the fourth necessary condition $x\mid (a+1)$ is immediately established because $[+a +1]$ is sandwiched between the two occurrences of $x$.  Furthermore, $\gcd(a,y)$ must divide the consecutive, hence relatively prime, integers $x-2$ and $x-1$, so $a$ and $y$ are relatively prime, establishing the last requirement.   

\emph{Part II, Sufficiency:}
Now suppose the five conditions of the theorem are met for some reduced set of four skip sizes.  As usual, we just need to establish the conditions of the master lemma for some proposed 7 step path.  We consider the pattern $[+a+1-x+1-y+1-x]$ which will be a cycle under the third condition, so it only remains to show that this pattern is realizable. Note that since $x$ and $y$ are each at least 3, the seven terms of this pattern must be distinct, so that ``weakly realizable'' implies ``realizable'' for that pattern.

The parity of the various intermediate paths is easily checked to be consonant with the requirements of the master lemma, once we have shown that $a, x,$ and $y$ are pairwise relatively prime.   The fourth condition $x \mid (a+1)$ shows that $x$ divides a number relatively prime to $a$, hence $\gcd(a,x) = 1$.    Moreover, the third condition can be rewritten to give that $(a + 1) - 2x = y-2$ showing that $x \mid (y-2)$.  As $y$ is odd, this shows $x$ divides a number relatively prime to $y$, hence $\gcd(x,y) = 1$.  The fifth condition is, of course, precisely the statement that the remaining pair is relatively prime, thus concluding our proof.\hfill$\Box$

\begin{claim}\label{OnlyOneEvenClass}
Given a set of skip sizes $\{a,b,x,y\}$ with $\gtwo ab\gtwo{}{}\eqtwo xy$, then this set produces an odd cycle if and only if it produces a 3-cycle.
\end{claim}

\noindent\textbf{Proof:}
Proof of sufficiency is trivial, so assume the skip set produces an odd cycle.
By Claim~\ref{fourskipsizes}, the cycle must be length 7 or less. By Claim~\ref{SevenCycle}, this skip set has the wrong pattern of equivalence classes to produce a 7-cycle.  If it produces a 3-cycle, we are done, so assume it produces a 5-cycle.

Claim~\ref{oddblocks} shows that skips $x$ and $y$ must occur in pairs in the 5-cycle, and thus can only occur once apiece by the No Pairs Corollary.  Because $\gtwo ab$ we see that $a$ does not divide $b$, hence the only possible pattern for the 5-cycle is $[babyx]$.  Applying the appropriate Master Lemma Corollaries, we see that this requires the signed pattern to have the form $[-b\pm a +b+y-x]$.  The Cycle Equation then necessitates that $\pm a = x-y$ and hence by Claim~\ref{skipsetsize3} the set $\{a,x,y\}$ produces a 3-cycle.
\hfill$\Box$

%
%
%

\begin{claim}\label{OddCycleFromTwoEvens} Given a reduced set of skip sizes $\{a, b, x, y\}$ with $\eqtwo{a}{b}\gtwo{}{}\eqtwo{x}{y}$ that does not produce a 3-cycle, then this set produces an odd cycle if and only if the following conditions are met under some appropriate labeling:
	\begin{itemize}
		\item  $b \mid a$; 
		\item  $\gcd(a,x) \mid b$;  
		\item $\gcd(a,y) \mid b$; and
		\item  $a=2b + (y-x)$.	
 	\end{itemize}
\end{claim}

\noindent\textbf{Proof:}  
\emph{Part I, Necessity:}
Suppose we have a skip set that produces an odd cycle but not a 3-cycle.  By Claim~\ref{SevenCycle}, we see that this set cannot produce a 7-cycle as it has the wrong equivalence class pattern.  As it cannot produce cycles longer than 7 (Claim~\ref{fourskipsizes}) it must produce a 5-cycle.  

As in the previous proof, we know that, $x$ and $y$ can only be used once each in a 5-cycle, and thus, without loss of generality, the cycle must be $[babyx]$.  The Divisibility condition applied to the subpaths $[bab]$,  $[aby]$, and $[xba]$ gives us that $b \mid a$, that $\gcd(a,y) \mid b$, and that $\gcd(a,x)\mid b$.  Making the cycle conform to the BPT requires the signs to alternate for the even numbers and again for the odd numbers, so that the cycle has pattern $[-b+a-b+y-x]$ or $[+b -a + b +y -x]$.  The  Cycle Equation then requires $a = 2b + (x-y)$ or $a = 2b + (y-x)$.  Interchanging the labels for $x$ and $y$ if necessary, we obtain the fourth condition.

\emph{Part II, Sufficiency:}
Suppose we have a (reduced) skip set $\{a,b,x,y\}$
where $x$ and $y$ are odd and $a$ and $b$ are in the same even equivalence class.  Suppose also that $b \mid a$; both $\gcd(a,x)$ and $\gcd(a,y)$ divide $b$; and $a=2b + y-x$.   We now consider the pattern $[+b -a + b+y-x]$.

The Divisibility condition for subpaths $[bab], [aby], [xba]$ all follow directly from these assumptions.  Further, $\gcd(b,x)$ and $\gcd(b,y)$ both divide $b-a$, hence $[+b-a+b+y]$ and $[-x+b-a+b]$ also satisfy the Divisibility condition.  By Claim~\ref{BothPaths}, the alternate subpaths $[yxb]$ and $[byx]$ also satisfy Divisibility.  Thus all subpaths of length 3 satisfy Divisibility.  

Because the signs alternate for the even skips, alternate for the odd skips, and the first odd skip is positive, the path $[+b-a+b+y-x]$ satisfies the BPT.  By Claim~\ref{SingleSeparatedNeighbors}, this means all subpaths of length 3 satisfy Parity as well.  Invoking Claim~\ref{BothPaths}, we see that all subpaths satisfy Divisibility and Parity.  So because the Cycle Equation is satisfied, we have a valid pattern for a 5-cycle. 
\hfill$\Box$

\begin{claim}\label{FiveCycleFromOneEven} Given a reduced skip set $\lbrace a, x, y, z \rbrace$ with $\gtwo ax\eqtwo{}{}\eqtwo yz$, and which produces no 3-cycles, then it produces a 5-cycle if and only if the following conditions are met for some appropriate labeling of $x,y,$ and $z$.
	\begin{itemize}
		\item  $a = \pm(2x-y-z)$; and 
		\item  $x\mid y$;  and
		\item  $\gcd(y,z) \mid x$; and
		\item  $\gcd(a,y) \mid x$.
 	\end{itemize}
\end{claim}

\noindent\textbf{Proof:}  
\emph{Part I: Necessity:}  Suppose this set produces a 5-cycle.  Because there have to be an even number of odd path elements (Claim~\ref{oddblocks}) exactly one of the odd numbers, say $x$, is used twice.  Taking into account relabeling and reversing of patterns, there are only two possible 5-cycles: $[axyzx]$ and $[axyxz]$.

To satisfy the BPT, the first pattern must be of the form $[\pm a + x-y+z-x]$.  The Cycle Equation then gives that $\pm a = y-z$, and thus the skip set produces a 3-cycle as well in accordance with Claim~\ref{skipsetsize3}.  Our conditions preclude this possibility, so we are only left with the 5-cycle $[axyxz]$.

Again applying the BPT, the cycle must have the form $[\pm a + x -y +x -z]$.  The Cycle Equation then gives us our first condition: $a = \pm(2x-y-z).$  The other three conditions are immediate consequences of the Divisibility conditions applied to various subpaths.

\emph{Part II: Sufficiency:}  
Assume $a = 2x-y-z$.  The other case, $a = -(2x-y-z)$ is handled similarly.  We will show that $[-a + x -y+x-z]$ is a valid cycle.  It already satisfies the Cycle Equation, so it only remains to show that all the Divisibility and Parity conditions are satisfied.

By Claims~\ref{BothPaths} and~\ref{SingleSeparatedNeighbors}, it suffices to check Divisibility for all length 3 subpaths and the BPT.  This technique was demonstrated in the proof of Claim~\ref{OddCycleFromTwoEvens}, and so the details will not be repeated here.  We note only that because $\gcd(a,x)$ divides $a,x,$ and $y$, the Cycle Equation shows that it also divides $z$.
\hfill$\Box$

We can now combine the forgoing results into the following theorem which completely characterizes all (reduced) sets of size 4 that produce odd cycles.

\begin{thm}[Skip Sets of Size 4]
A reduced set of four numbers forces discrepancy two if and only if at least one of the following is true:
\begin{itemize}
\item There is an even number $a$ and odd numbers $b$ and $c$ such that $a+b = c$; or
\item There are two even numbers $a$ and $b$, where $b$ divides $a$ and the odd numbers $x$ and $y$ are relatively prime to $a$, with 
$a = 2b +y-x$; or
\item The set is $\{a,x,y,z\}$ where $a$ is the only even number, $x\mid y$, $y$ is relatively prime to $a$ and $z$, and they all satisfy the equation \\ $a = \pm(2x-y-z)$; or
\item The set is $\{a,x,y,1\}$ where $a$ is the only even number, $x$ divides $a+1$, $y$ is relatively prime to $a$, and they all satisfy the equation $a=2x+y-3$.
\end{itemize}
\end{thm}

\noindent\textbf{Proof:} 
Sufficiency follows immediately from  Claims~\ref{skipsetsize3},~\ref{OddCycleFromTwoEvens},
~\ref{FiveCycleFromOneEven}, and~\ref{SevenCycle}.

To see necessity, we begin by noting that Claim~\ref{fourskipsizes} shows we cannot get odd cycles longer than 7.

Claim~\ref{skipsetsize3} shows that the first condition is the only way to get 3-cycles.  Similarly, Claim~\ref{SevenCycle} shows the last condition is the only way to get 7-cycles.  So assume that our skip set only produces a 5-cycle.

Claim~\ref{lowestOddClass} shows that the skip set must have at least two odd numbers in it in order to produce a cycle.  

If there is only one even number in the skip set, then Claim~\ref{FiveCycleFromOneEven} shows that the third condition is the only way this can happen.

Finally, if the set has two even numbers in it, then they must be in the same equivalence class by Claim~\ref{OnlyOneEvenClass}, else it would produce a 3-cycle.  But then Claim~\ref{FiveCycleFromOneEven} completes our argument.
\hfill$\Box$

\section{Complexity}\label{Complexity}

The characterization of skip sets of size 3 that force discrepancy two was relatively straightforward, whereas that for skip sets of size 4 was {\em much} messier.  Some computational investigations suggest that path and odd cycle lengths grow exponentially with skip set size, as shown in Figure~\ref{table}, leading us to believe that getting a handle on odd cycles will continue to grow rapidly in difficulty. We therefore define the following problem:\\

\begin{figure}
\noindent\begin{tabularx}{\linewidth}{ |c|c|c|c|X| }
\hline
$|S|$ & Structure & len & Start & Example\\ \hline
2 & cycle &- &-&-\\
2 &	path	 & 3 & 1 & [1 3 1]\\ \hline
3 & cycle & 3 & 0 & [2 1 3]\\
3 & path & 7 &  11 & [1 5 1 7 1 5 1]\\ \hline
4 & cycle & 7 & 48 & [8 1 3 1 5 1 3]\\ 
4 & path & 18 & 70 & [7 1 4 9 1 4 7 1 9 1 7 1 4 9 1 4 7 1]\\ \hline
5 & cycle & 19 & 45756 & [62 1 9 5 11 1 5 1 9 1 11 1 9 1 5 1 11 5 9]\\ 
5 & path & 53 &  15962 & [2 35 1 2 1 9 37 1 9 1 2 35 9 2 37 1 2 9 1 35 1 2 37 1 2 35 1 2 37 1 2 35 1 2 9 1 37 9 2 35 1 9 1 2 37 1 2 1 9 35 1 2 1] \\ \hline
6 & cycle & 47 & 4836 &[6 3 1 4 121 1 4 1 3 6 4 1 13 1 3 6 4 1 3 1 13 6 4 1 3 6 4 1 13 4 1 3 6 4 1 3 6 3 13 1 3 4 1 3 6 1 13] \\
6 & path & 165 & 2848  &  [4 1 3 10 1 3 6 59 1 10 6 3 1 4 6 3 1 10 3 1 4 6 3 1 4 10 1 3 6 3 1 4 10 59 1 4 6 3 1 10 3 1 4 6 3 1 4 10 1 3 6 3 1 4 6 3 1 59 1 3 1 4 6 3 1 10 3 1 4 6 3 1 4 10 1 3 6 3 1 4 6 59 1 4 6 3 1 4 6 3 1 10 3 1 4 6 3 1 4 10 1 3 6 3 1 4 6 3 1 4 59 3 6 3 1 4 6 3 1 10 3 1 4 6 3 1 4 10 1 3 6 3 1 4 3 1 59 1 6 3 1 4 6 3 1 10 3 1 4 6 3 1 4 10 1 3 6 3 1 4 6 3 1 4 1]\\
\hline
\end{tabularx}\caption{Long paths and cycles for various sizes of skip set $S$. The signs in the patterns can be inferred as $+$ or $-$ depending on whether the term is an even or odd multiple resp. of the next skip.}\label{table}
\end{figure}

\noindent{\sc Discrepancy One (D1)}\\
{\sc Instance:} A set $S=\{s_1, s_2, \ldots, s_n\}$ of $n$ positive integers.\\
{\sc Question:} Can the natural numbers be 2-colored so that homogeneous arithmetic progressions with differences in $S$ all have discrepancy 1?\\

We don't know whether D1 is in co-NP, because cycle sizes seem to grow exponentially with the size of S. And we don't see how to show D1 is in NP, because showing a positive answer seems to require constructing a 2-coloring of the interval $[0, 2\cdot\mbox{lcm}(S)-1]$, which may have length, again, exponential in the size of $S$. But we can show that D1 is NP-hard.

\begin{thm}\label{complexity}
{\sc Discrepancy One} is NP-hard.
\end{thm}

\noindent\textbf{Proof:} We do this by reduction from the following problem, shown to be NP-complete by Cieliebak et al. in \cite{equalsum}:\\

\noindent{\sc Equal Sum Subsets of Different by One Cardinality (ESS)}\\
{\sc Instance:} A set $A=\{a_1, a_2, \ldots, a_n\}$ of $n$ positive integers.\\
{\sc Question:} Are there two disjoint subsets $X, Y\subseteq A$ with
$|X| = |Y|+1$ such that the sum of the elements of $X$ equals the sum
of the elements of $Y$?\\

Let $A=\{a_1, a_2, \ldots, a_n\}$ be an instance of ESS with elements given in increasing order. We transform this to an instance of D1 by creating the skip set $S=\{s_1, s_2\ldots, s_n, M, t\}$ where 
\begin{equation}\label{transformation}
\begin{split}
M   &= n\prod\limits_{1\leq i < j \leq n} (a_j-a_i)
\prod\limits_{i=1}^n (na_i+1)\\
r &= \mbox{least multiple of $n$ greater than } (n/2)(n(a_n-a_1) + 1)+1\\
s_i &= nMa_i+rM+1\\
t   &= (r-1)M+1.
\end{split}
\end{equation}
We now show that ESS has a positive answer if and only if the corresponding instance of D1 has a negative answer.

We first suppose that ESS has a positive answer. This means that we can find disjoint subsets $X, Y\subseteq A$ with $|X| = |Y|+1$ such that the sum of the elements of $X = \{a_{i_1},\ldots a_{i_{k+1}}\}$ equals the sum of the elements of $Y=\{a_{j_1}, \ldots, a_{j_k}\}$ for some $k$. We show that the following pattern yields an odd cycle in the corresponding skip graph: 
$[+s_{i_1}-s_{j_1} + s_{i_2}-s_{j_2} + \ldots + s_{i_k}-s_{j_k}+ s_{i_{k+1}}-t-M]$. We first note that it satisfies the cycle equation, since
\begin{equation}\label{cyclesum}
\begin{split}
\mbox{sum}& = +s_{i_1}-s_{j_1} + s_{i_2}-s_{j_2} + \ldots + s_{i_k}-s_{j_k}+
s_{i_{k+1}}-t-M\\
& = nM\cdot[(a_{i_1}+\ldots + a_{i_{k+1}}) - 
(a_{j_1}+ \ldots+ a_{j_k})]\\
& + (k+1)(rM+1) - k(rM+1) - ((r-1)M + 1) - M\\
& = nM\cdot 0 + (rM+1) - (r-1)M - 1 - M\\
& = 0
\end{split}
\end{equation}

Next we show that the elements of $S$ are pairwise relatively prime so that the Divisibility condition is satisfied.\\
{\em Case 1:} $\gcd(s_i, s_j)$. Suppose $d\mid s_i$ and $d\mid s_j$. Then $d\mid (s_i-s_j) = (nMa_i+rM+1)-(nMa_j+rM+1) = nM(a_i-a_j)$. By the definition of $M$, any prime dividing $n$ or $(a_i-a_j)$ also divides $M$, so any prime dividing $d$ must also divide $M$, and thus cannot divide $s_i=Ma_i+rM+1$, so $d=1$.\\
{\em Case 2:} $\gcd(s_i, M)=1$ since $s_i$ is one more than a multiple of $M$.\\
{\em Case 3:} $\gcd(M, t)=1$ since $t$ is one more than a multiple of $M$.\\
{\em Case 4:} $\gcd(s_i, t)$. Suppose $d\mid s_i$ and $d\mid t$. Then $d\mid (s_i-t) = (nMa_i+rM+1)-((r-1)M+1) = M(na_i+1)$. And again, any prime dividing $M(na_i+1)$ must also divide $M$, and therefore cannot divide $s_i$, so $d=1$.

It is straightforward to verify the Parity condition for the pattern, so that the Master Lemma applies, and the pattern is weakly realizable. This odd cycle may repeat arcs or terms, but it is sufficient to prove that no 2-coloring exists in the skip graph, so that our instance of D1 has a negative answer.

Now we prove the other direction. Suppose that our instance of D1 has a negative answer. Then there is some odd cycle $C$ in the skip graph induced by the set $S=\{s_1, \ldots, s_n, M, t\}$. All members of that set are odd except for $M$, so that $C$ must have the form 
$[(\pm M) + x_1 - x_2\ (\pm M) + x_3 - x_4\ (\pm M) + \cdots + x_{2l-1} - x_{2l}]$, where each $x_i$ is an (odd) element of $S-\{M\}$, and each parenthesized $(\pm M)$ expression indicates that the term may or may not be present in the pattern, and if present may be $+M$ or $-M$. But since the cycle has odd length, there need to be an odd number of $\pm M$s in the pattern. \\
{\em Claim:} If $i\neq j$, then $x_i \neq x_j$. That is, each skip size other than $M$ can occur at most once in the cycle.\\
{\em Proof:} Suppose some skip size occurs twice. Then select $i$ and $j$ so that $i<j$ and $j-i$ is minimized over all pairs where $x_i = x_j$. The path from the left endpoint of the $x_i$ arc to the left endpoint of the $x_j$ arc must have the form $[(\pm M) + x_{i_1} - x_{i_2}\ (\pm M) + x_{i_3} - x_{i_4}\ (\pm M) + \cdots + x_{i_{2m-1}} - x_{i_{2m}}\ (\pm M)]$, with no repeated arcs other than, perhaps, $M$. Let $\Delta = a_n - a_1$. Since 
$|x_i-x_j| = nM|a_i-a_j| \leq nM\Delta$, and $M$ appears at most $n/2$ times, the distance between the start and end of that path is at most 
$(n/2)nM\Delta + (n/2)M = (nM/2)(n\Delta + 1)$. 
$r$ was selected so that this distance is less than $2x_i$, giving a contradiction, since the left endpoints of two different $x_i$-arcs must be at least $2x_i$ apart.\\
{\em Claim:} The cycle contains exactly one $M$-arc and exactly one $t$-arc.\\
{\em Proof:} We consider the skip sizes of the cycle modulo $nM$. We have

\begin{center}
\begin{tabular}{lll}
$s_i$ &$\equiv 1$    &$\pmod{nM$}\\
$t$   &$\equiv -M+1$ &$\pmod{nM}$\\
$M$   &$\equiv M$    &$\pmod{nM}$
\end{tabular}
\end{center}

The $M$-arc must occur at least once since the cycle has an odd number of arcs, and can occur at most $n/2$ times. The $t$-arc can occur at most once, and there can be at most $n$ $s_i$-arcs. There is therefore no way for the values shown on the right sides of the congruences to add to any multiple of $nM$ other than 0. There must therefore be exactly one $M$-arc, exactly one $t$-arc, and exactly one more $s_i$-arc with a negative sign than with a positive sign. So our cycle must have the following distinct arcs in some order:
\begin{center}
$+s_{i_1} +s_{i_2} +\cdots +s_{i_l}$\\
$-s_{j_1} -s_{j_2} -\cdots -s_{j_l} -s_{j_{l+1}}$\\
$M$\\
$t$\\
\end{center}

Since the sum of these terms must be 0, we have, just as in equation (\ref{cyclesum}) 
$(a_{i_1}+\ldots + a_{i_{l}}) - 
(a_{j_1}+ \ldots+ a_{j_l} + a_{j_{l+1}}) = 0$.
showing that our instance of ESS has a positive answer.

Finally, we note that translating ESS to D1 requires performing the computations shown in the equations in (\ref{transformation}), which are clearly of polynomial time and space.\hfill$\Box$

\bibliography{DiscrepancyOne}{}
\bibliographystyle{plain}

\end{document}